\documentclass[titlepage,10pt]{article}
\usepackage{a4,epsf,amscd,amsfonts,amsmath,amssymb,amstext,amsthm,
latexsym,times}

\begin{document}

\renewcommand{\theenumi}{\roman{enumi}}   
\renewcommand{\thefootnote}{\fnsymbol{footnote}} 
\newcommand{\gc}{G^\mathbb C}                    
\newcommand{\qd}{\hfill{\( \square \)}}
\newcommand{\vn}{\mathbb C^{(1+n) \times N}}           
\newcommand{\mo}{ / \negmedspace /} 
\newcommand{\lo}{\mathrm{SO}_\mathbb R(1,n)^0}
\newcommand{\loc}{\mathrm{SO}_\mathbb C(1,n)}
\newcommand{\RR}{\mathbb R}
\newcommand{\CC}{\mathbb C}
\newcommand{\kn}{\mu^{-1}(0)}  
\newtheorem{theo}{Theorem}[section]              
\newtheorem{prop}{Proposition}[section]      
\newtheorem{lemm}{Lemma}[section]                  
\newtheorem{corr}{Corollary}[section]          
\newtheorem{defi}{Definition}[section]          
\newtheorem{beme}{Remark}[section]                
\newtheorem{beis}{Example}[section]              
\newtheorem{fact}{Fact}

\vskip5cm

\Large{\centerline{The extended future tube conjecture for SO(1,n)}}
\normalsize
\bigskip
\bigskip

\bigskip 
\centerline{P\footnotesize{ETER} \normalsize H\footnotesize{EINZNER 
 AND } \normalsize P\footnotesize{ATRICK} \normalsize 
S\footnotesize{CH\"UTZDELLER} \footnote{Supported by the SFB 237 of the DFG and
 the Ruth und Gerd Massenberg Stiftung } }

\bigskip
\bigskip
\normalsize

\newenvironment{Quote}%
 {\begin{list}{}{%
   \setlength{\rightmargin}{\leftmargin}}%
  \item[] \makebox[0pt][r]{}\ignorespaces}
 {\unskip\makebox[0pt][l]{}\end{list}}

\vskip1cm
\begin{quote}
A\footnotesize BSTRACT. \normalsize \
Let $C$ be the open upper light cone in $\mathbb R^{1+n}$ with respect to 
the Lorentz product. The connected linear Lorentz group 
$\lo$ acts on $C$ and therefore diagonally on the 
$N$-fold product $T^N$ where $T = \RR^{1+n} + iC \subset \CC^{1+n}.$
We prove that the extended future tube 
$\loc \cdot T^N$ is a domain of holomorphy.

\end{quote}

\vskip1cm

\noindent
For  $\mathbb K \in \{ \RR, \CC\}$ let $ \mathbb K^{1+n}$ denote
 the $(1+n)$-dimensional 
Minkowski space, i.e., on  $\mathbb K^{1+n}$ we have given 
the 
bilinear form $$(x,y) \mapsto x \bullet y := x_0y_0 -x_1y_1 - \cdots - x_ny_n$$
 where $x_j$ respectively $y_j$ are the components  of $x$ respectively $y$ in
 $\mathbb K^{1+n}.$ 
The group $\mathrm O_\mathbb K(1,n) = \{ g \in \mathrm{Gl}_\mathbb K(1+n) ;
gx \bullet gy = x \bullet y$ for all $x,y \in \mathbb K^{1+n}\}$ is called the
 linear Lorentz group. For $n \geq 2$ the group $\mathrm O_\RR(1,n)$
 has four connected components
 and $\mathrm O_\CC(1,n)$ has two connected components. 
The connected component of the
 identity $\mathrm O_\mathbb K(1,n)^0$ of $\mathrm O_\mathbb K(1,n)$ will be
  called the connected linear 
 Lorentz group. Note that $\mathrm {SO}_\RR(1,n) = \{ g \in 
\mathrm O_\RR(1,n);  \det(g) =1\}$ has two connected components and 
$\mathrm O_\RR(1,n)^0 = \lo$. In the complex case we have 
$\loc = \mathrm O_\CC(1,n)^0$. 

\medskip
\noindent
The forward cone $C$ is by definition the set 
$C := \{ y \in \RR^{1+n} ; y \bullet y > 0$ and $ y_0 > 0\}$ and
the future tube $T$ is the tube domain over $C$ in $\CC^{1+n}$, i.e., $T = 
\RR^{1+n} + iC \subset \CC^{1+n}$. Note that $T^N = T \times \cdots \times
T$ is the tube domain in the space  of complex $(1+n) \times N$-matrices 
$\CC^{(1+n) \times N}$ over $C^N = C \times \cdots \times C \subset
\RR^{(1+n) \times N}$. The group $\loc$ acts by matrix multiplication on 
$\vn$ and the subgroup $\lo$ stabilizes $T^N$. In this note we prove the 

\medskip
\noindent
\underline{Extended future tube conjecture:}
$$\loc \cdot T^N = \bigcup_{g \ \in \  \loc} g \cdot T^N \mbox{ is a domain of 
holomorphy}. $$ 

\medskip
\noindent
This conjecture arise in the theory of quantized fields for about 50 years.
We refer the interested reader to the literature
(\cite{HW}, \cite{Jo}, \cite{SV}, \cite{SW}, \cite{W}).
There is a proof of this conjecture in the case where $n=3$ (\cite{H2}),
\cite{Z}). The proof there uses essentially that $T$ can be realized as the 
set $\{Z \in \CC^{2\times 2} ; \frac 1{2i} ( Z - ^t \! \! \! \bar Z) $ 
is positive 
definite$\}$. Moreover the proof for $n=3$ is unsatisfactory. It does not 
give much information about $\loc \cdot T^N$ except for holomorphic convexity.

\medskip
\noindent
Here we prove that more is true. Roughly speaking, we show that the basic
 Geometric Invariant Theory results known for compact groups 
(see  \cite{H1}) also holds for $X:=T^N$ and the non compact group 
$\lo$. More precisely this means $\loc \cdot X = Z$ is a universal 
complexification of the $G$-space $X$, $G = \lo$, in the sense
of \cite{H1}. There exists complex analytic quotients  $X \mo G$ and 
$Z \mo \gc$, $\gc = \loc$, given by the algebra of invariant holomorphic 
functions and there is a $G$-invariant strictly plurisubharmonic 
function  $\rho:  X \to \RR$, which is an exhaustion on $X/G$. Let 
$$\mu : X \to \mathfrak g^\ast, \quad  \mu(z) (\xi) = 
\frac {d} {dt} \Big \vert_{t=0} ( t \to \rho(\exp it\xi \cdot z)),$$
be the corresponding moment map. Then the diagram

$$ \begin{matrix} \mu^{-1}(0) & \hookrightarrow & \quad X &  \quad 
 \hookrightarrow & \ \   Z \cr
                      & & & &  \cr
       \downarrow &  & \qquad \downarrow \pi & & \qquad \downarrow \pi^\CC \cr
                      & & & &  \cr
                   \mu^{-1}(0)/G  & \equiv & \quad   X \mo G& \quad 
 \equiv & \quad \ \  Z \mo \gc
\end{matrix}$$

\medskip
\noindent
where all maps are induced by inclusion is commutative, 
$X / \!  /   G, X, Z$ and 
$Z \mo \gc$ are Stein spaces and $\rho \vert \kn$ induces a 
strictly plurisubharmonic exhaustion on $\kn/G =$ $X / \! \negmedspace / G$
 $=$ $Z \mo \gc.$ Moreover the same statement holds if we replace $X=T^N$
with a closed $G$-stable analytic subset $A$ of $X$.

\section{Geometric Invariant Theory of Stein spaces} \label{GIT}

\medskip           
\noindent
Let $Z$ be a Stein space and $G$ a real Lie group acting as a group of 
holomorphic transformations  on 
$Z$. A complex space $Z \mo G$ is said to be an analytic Hilbert quotient
of $Z$ by the given $G$-action if there is a $G$-invariant surjective 
holomorphic map $\pi : Z \to Z \mo G$, such that for every open Stein 
 subspace  $Q \subset  Z \mo G$

\begin{enumerate}

\item its inverse image $\pi^{-1}(Q)$ is an open  Stein 
subspace of $Z$ and

\item
$ \pi^\ast \mathcal O_{Z \mo G} (Q) = \mathcal O(\pi^{-1}(Q))^G$, where
$\mathcal O(\pi^{-1}(Q))^G$ denotes the algebra of $G$-in\-variant holomorphic
 functions on $\pi^{-1}(Q)$ and $\pi^\ast$ is the pull back map. 

\end{enumerate}

\noindent Now let $G^c$ be a linearly reductive complex Lie group.
A complex space  $Z$ endowed with a holomorphic action of $G^c$ is called a 
holomorphic $G^c$-space.

\begin{theo} \label{hilbert quotient}
Let $Z$ be a holomorphic $G^c$-space, where $G^c$ is a linearly reductive 
complex Lie group.
\begin{enumerate}
\item If $Z$ is a Stein space, then the analytic Hilbert quotient 
$Z \mo G^c$ exists and is a Stein space.

\item If $Z \mo G^c$ exists and is a Stein space, then $Z$ is a Stein space.
\end{enumerate}
\end{theo}

\medskip
\noindent
\textit{Proof.} Part i. is proven in \cite{H1} and part ii. in \cite{HMP}.
\qd

\begin{beme} \ 

\begin{enumerate}
\item If the analytic Hilbert quotient $ \pi : Z \to Z \mo G^c$ exists, then 
every fiber $\pi^{-1}(q)$ of $\pi$ contains a unique $G^c$-orbit $E_q$
of minimal dimension. Moreover, $E_q$ is closed and $\pi^{-1}(q) = \{ z \in Z;
 E_q \subset \overline{G^c . z} \}$. Here $ \ \bar \ \bar  \  \   $ denotes the topological 
closure. 

\item Let X be a subset of $Z$, such that 
$G^c \cdot X := \bigcup_{g \in G^c} g \cdot X = Z$ and assume that
$Z \mo G^c$ exists. Then $G^c \cdot X$ is a Stein space if and only if 
$Z \mo G^c = \pi(X)$ is a Stein space.

\item Let $V^c$ be a finite dimensional complex vector space with a 
holomorphic linear action of $G^c$. Then the algebra $\CC[V^c]^{G^c} $ of 
invariant polynomials is finitely generated (see e.g. \cite{Kr}).
\end{enumerate}
\end{beme}

\medskip
\noindent
In particular, the inclusion $ \CC[V^c]^{G^c} \hookrightarrow \CC[V^c]$
defines an affine variety $ V^c \mo G^c $
 and an affine morphism $\pi^c: V^c \to
 V^c \mo G^c$. If we regard $V^c \mo G^c $ as a complex space, then $ \pi^c : 
V^c \to V^c \mo G^c$ gives the analytic Hilbert quotient of $V^c$ 
(see e.g. \cite{H1}).

\begin{beme}
For a non-connected linearly reductive complex
group $G$ let $G^0$ denote the connected 
component of the identity and let $Z$ be a holomorphic $G$-space. The analytic
Hilbert quotient $Z \mo G$  exists if and only if the quotient $Z \mo G^0$
exists. Moreover, the quotient 
map $\pi_G : Z \to Z \mo G$ induces a map 
$\pi_{G/G^0} : Z \mo G^0 \to Z \mo G$ which is finite.
In fact the diagram
$$ \begin{matrix} 
         &                      & Z               &                 &       \cr
         &  \pi_{G^0} \swarrow  &                 & \searrow  \pi_G &       \cr
         &                      &                 &                 &       \cr
Z\mo G^0 &                      & \longrightarrow &                 &Z\mo G \cr
         &                      & \pi_{G/G^0}     &                 &
\end{matrix}$$
commutes and $\pi_{G/G^0}$ is the quotient map for the induced action of
the finite group $G/G^0$ on $Z \mo G^0$. 
\end{beme}

\section{The geometry of the Minkowski space} \label{Minkowski}

\medskip
\noindent
Let $\mathbb K$ denote either the field $\RR$ or $\CC$ and 
$(e_0,..,e_n)$ the standard orthonormal basis for $\mathbb K^{1+n}$.
The space $\mathbb K^{1+n}$ together with
the quadratic form $\eta(z) = z_0^2-z_1^2 - \cdots -z_n^2$, where $z_j$
are the components of $z$, is called the $(1+n)$-dimensional linear Minkowski
 space.  Let $<,>_L$ denote the 
symmetric non-degenerated bilinear form which corresponds to $\eta$, i.e.,
 $ z \bullet w :=<z,w>_L =$ $ ^t \! zJw$ where 
$ ^t \! z$ denotes the transpose of $z$ and
$J = (e_0,-e_1,...,-e_n)$ or equivalently $z \bullet w = <z,Jw>_E$ 
where $<,>_E$
denotes the standard Euclidean product on $\RR^{1+n}$, respectively its 
$\CC$-linear extension to $\CC^{1+n}$.

\medskip
\noindent
Let $\mathrm O_\mathbb K(1,n)$ denote the subgroup
of Gl$_\mathbb K(1+n)$ which leave $\eta$ fixed, i.e., 
$\mathrm O_\mathbb K(1,n) = \{ g \in \mathrm{Gl}_\mathbb K(1+n) ; 
g z \bullet gw = z \bullet w$ for all $ z,w \in  \mathbb K^{1+n}\}$. 
Note that $\mathrm {SO}_\mathbb K(1,n) = \{ g \in \mathrm O_\mathbb K(1,n) ;
\det g= 1\}$ is an open subgroup of $\mathrm O_\mathbb K(1,n)$. 
For $\mathbb K = \CC$,
$\loc$ is connected. But in the real case $\mathrm {SO}_\mathbb R(1,n)$ 
consists of two 
connected components $(n \geq 2)$. The connected component 
$\lo = \mathrm O_\mathbb R(1,n)^0$ of the identity  is called the connected 
linear Lorentz group. Note that $\lo$ is not an algebraic subgroup of 
$\mathrm {SO}_\mathbb R(1,n)$ but is Zariski dense in 
$\mathrm {SO}_\mathbb R(1,n)$.
We have 
$\mathbb K[\eta]= \mathbb K [ \mathbb K^{1+n}]^{\mathrm {SO}_\mathbb K(1,n)} =
\mathbb K[ \mathbb K^{1+n}]^{\mathrm {O}_\mathbb K(1,n)}$.

\medskip
\noindent
Now let $\mathbb C^{(1+n) \times N} = \mathbb C^{1+n} \times \cdots \times 
 \mathbb C^{1+n}$ be the $N$-fold product of  $\mathbb C^{1+n}$, i.e., the 
space of complex $(1+n) \times N$- matrices.
The group $\mathrm O_\mathbb C(1,n)$ acts on 
$\mathbb C^{(1+n) \times N}$ by left 
multiplication. A classical result in Invariant Theory says that 
$\mathbb C[ \mathbb C^{(1+n) \times N}]^{\mathrm {O}_\mathbb C(1,n)}$ 
is generated by the polynomials $p_{kj}(z_1,..,z_N) = z_k \bullet z_j$
where  $z =(z_1,..,z_N) \in  \mathbb C^{(1+n) \times N}$.

\begin{beme}
The (algebraic) Hilbert quotient $ \mathbb C^{(1+n) \times N} \mo 
\mathrm {O}_\mathbb C(1,n)$ can be identified with the space 
$\mathrm{Sym_N}(\mathrm{min} \{ 1+n,N \})$ of symmetric 
$N \times N$-matrices of 
rank smaller or equal $\mathrm{min} \{ 1+n,N \}$.
\end{beme}

\medskip
\noindent
With this identification the quotient map $\pi_\CC :
  \mathbb C^{(1+n) \times N}
\to  \mathbb C^{(1+n) \times N} \mo 
\mathrm {O}_\mathbb C(1,n)$ is given by $\pi_\CC (Z) =  {}^t \! ZJZ$
where $ ^t \! Z$ denotes the transpose of $Z$ and $J$ is as above.
For the group $\mathrm {SO}_\mathbb C(1,n)$ the situation is slightly more 
complicated. If $N \geq 1+n$ additional invariants appear, but they are not 
relevant for our considerations, since the induced map 
$  \mathbb C^{(1+n) \times N} \mo \mathrm {SO}_\mathbb C(1,n) \to 
\mathbb C^{(1+n) \times N} \mo \mathrm {O}_\mathbb C(1,n)$ is finite.

\medskip
\noindent
There is a well known characterization of closed $\mathrm O_\CC(1,n)$-orbits 
in 
$ \mathbb C^{(1+n) \times N}.$ In order to formulate  this we need more
notations. Let $z=(z_1,..,z_N) \in \vn$ and   
$L(z) := \CC z_1 + \cdots +  \CC z_N $ be the  subspace 
of $\mathbb C^{1+n}$ spanned by $z_1,..,z_N$. The Lorentz  product
$<,>_L$ restricted to $L(z)$ is in general degenerated.
Thus let  $L(z)^0= \{ w \in L(z); <w ,v>_L = 0$ for all $ v \in L(z) \}$.
It follows that dim$L(z)/L(z)^0 =$ rank$  ( ^t \!zJz) =$
rank $\pi_\CC(z)$. Elementary consideration show the following.

\begin{lemm} \label{closedorbits}
The orbit $\mathrm O_\CC(1,n) \cdot z $ through $z \in \vn$ is closed
if and only if the orbit  $\loc \cdot z$ is closed and this 
is the case if and only if
$L(z)^0 = \{ 0 \} $, i.e., dim$L(z)$ = rank $\pi_\CC(z)$.  \qd
\end{lemm}
 
\medskip
\noindent
The light cone  $ N:= \{ y \in \mathbb R^{1+n} ; 
\eta(y) = 0\}$ is of codimension one and its complement $\RR^{1+n} \backslash
N$ consists of three connected components 
(here of course we assume $n \geq 2$). By 
the forward cone $C$ we mean the connected component 
which contains $e_0$. It is easy to see that
$C = \{ y \in \mathbb R^{1+n}; y \bullet e_0 > 0 \mbox{ and }\eta(y) > 0\}$
$ =\{ y \in \mathbb R^{1+n}; y \bullet x > 0\mbox{ for all } x \in N^+\}$
where $N^+ = \{ x \in N; x \bullet e_0  > 0\}. $ In particular, $C$ is 
an open convex cone in $\RR^{1+n}$. Since $J$ has only one positive 
Eigenvalue, the following version of the Cauchy-Schwarz inequality holds.

\begin{lemm} \label{CW}
If $\eta(y) > 0$, then $\tilde x \bullet y \leq 0$ for  $\tilde x := x- \frac 
{x \bullet y} { \eta (y)^2} y $ and all $x \in \RR^{1+n}$. In particular
$$\eta( x) \cdot  \eta( y) \leq (x \bullet y)^2$$
 and equality holds if and only if $x$ and $y$ are linearly dependent.
\qd
\end{lemm} 

\medskip
\noindent
The elementary Lemma has several consequences which are used later on. 
For example,

\begin{itemize}

\item if $y_1, y_2 \in C^\pm := C \cup (-C) = \{y \in \RR^{1+n} ; \eta(y)>0\}$,
then $y_1 \bullet y_2 \neq 0$. Moreover,

\item if $y_1,y_2 \in N = \{ y \in \RR^{1+n}; 
\eta(y) =0\}$, and $y_1 \bullet y_2 =0$, then 
$y_1$ and $y_2$ are linearly  dependent.

\end{itemize}

\medskip
\noindent
The tube domain $T = \mathbb R^{1+n} + iC \subset \mathbb C^{1+n}$ over $C$ is
 called the future tube. Note that $\lo$ acts on $T$ by $g \cdot (x+iy) =
gx +i gy$ and therefore on the $N$-fold product $T^N = T \times \cdots \times
 T \subset \vn$ by matrix multiplication.

\begin{beme}
It is easy to show that the $\lo$-action on $C$  and consequently also on 
$T^N$ is proper. In particular $T^N/\lo$ is a Hausdorff space. 
\end{beme}

\medskip
\noindent
The complexified group $\loc$ does not stabilize $T^N$. The domain 
$$\loc \cdot T^N = \bigcup_{ g \in \loc} g \cdot 
 T^N$$ is called the extended future 
tube.

\section{Orbit connectedness of the future tube} \label{oconnected}

\medskip           
\noindent
Let $G$ be a Lie group acting on  $Z$. A subset 
$X \subset Z$ is called  orbit connected with respect to the $G$-action 
on $Z$ if $\Sigma(z) = \{ g \in G ;
g \cdot z \in X \}$ is connected for all $z \in X$. 
\\
In this section we prove the following

\begin{theo} \label{orbit-connected}
The $N$-fold product $T^N$ of the  future tube 
 is orbit connected with respect to 
the $\loc$-action on $\vn$.

\end{theo}

\medskip
\noindent 
We first reduce the proof of this Theorem for the $\loc$-action
to the proof of the  related statement about the Cartan subgroups of $\loc$. 
For this we use the results of Bremigan in \cite{Brem}. 
For the convenience of the reader we briefly recall those parts, which are 
relevant for the proof of Theorem~\ref{orbit-connected}. 

\medskip
\noindent
Starting with a simply connected complex semisimple Lie group $G^\mathbb C$ 
with a given real form $G$ defined by an anti-holomorphic group involution, $g 
\mapsto \bar g$, there is a subset $S$ of $G^\mathbb C$ such that 
$GSG$ contains an open $G \times G$-invariant dense subset of $ G^\mathbb C$. 
The set $S$ is given as follows.

\medskip
\noindent
Let $Car(G^\mathbb C) = \{H_1,..,H_\ell\}$ be a complete set of representatives
 of the Cartan subgroups of $G^\mathbb C$, which are defined over
 $\mathbb R$.
 Associated to each $ H \in Car(G^\mathbb C$) are the Weyl group 
$\mathcal W(H) := N_{G^\mathbb C}(H)/H$, the real Weyl group 
$ \mathcal W_\mathbb R(H) := \{ gH \in \mathcal W(H)
; \bar gH = gH \}$ and the totally real Weyl group 
$ \mathcal W_{\mathbb R!}(H) := 
\{ gH \in \mathcal W_\mathbb R(H) ; \bar g =g\}$. Here $N_{\gc}(H)$ denotes the
normalizer of $H$ in $\gc$.

\medskip
\noindent
For $H \in Car(G^\mathbb C)$  let $R(H)$ be a 
complete set of representatives of the
 double coset space   
$\mathcal W_{\mathbb R!}(H)\backslash \mathcal W_\mathbb R(H)/
\mathcal W_{\mathbb R!}(H)$
chosen in such a way that $\bar \epsilon = \epsilon^{-1}$ holds for all 
$\epsilon \in Car(G^\mathbb C)$. Then $S := \cup H  \epsilon  $ has the 
claimed properties.

\medskip
\noindent
Although $\mathrm {SO}_\mathbb C(1,n)$ 
is not simply connected, the results above
remain true for  $G:= \mathrm {SO}_\mathbb R(1,n)^0$ and 
$G^\mathbb C := 
\loc$, as one can see by going over to the universal covering. 

\begin{beme} Using the 
classification of the $\lo \times \lo$-orbits in $\loc$ as 
presented in \cite{Jo}, the same result can be obtained for $\gc = \loc$.
\end{beme}

\medskip
\noindent
Since $T^N$ is $\lo$-stable, $\lo$ is connected and $\lo \cdot  S \cdot \lo$ 
is dense in $\loc$, 
Theorem~\ref{orbit-connected} follows
from 

\begin{prop}
The set $\Sigma_S(w) := \{ g \in S ; g \cdot w \in T^N\}$ is connected for
all $w \in T^N$. 
\end{prop}

\medskip
\noindent
In the case $n = 2m - 1$ we may choose $Car(\loc)
 = \{H_0\} $
 where
$$H_0 = \left\{ \left( \begin{matrix} 
\sigma &  0 & \cdots & 0 \cr
     0      & \tau_1 & \ddots & \vdots \cr
 \vdots & \ddots & \ddots & 0 \cr
       0    & \cdots & 0 &   \tau_{m-1}   \end{matrix} \right) ;  \sigma \in 
\mathrm {SO}_\mathbb C(1,1), \tau_j \in \mathrm {SO}_\mathbb C(2)
\right\} \mbox{  and
 } R(H_0) = \{\mathrm{Id}\}.$$
In the even case $n = 2m$ we make the choice $Car(\loc)
 = \{H_1,H_2\} $ where
$$ H_1 = \left \{ \left (  \begin{matrix} 
    h        &       0         \cr
    0       &        1      \end{matrix} \right) 
; h \in H_0 \right\},
 H_2 = \left \{ \left (  \begin{matrix} 
   1        &     0  & \cdots & 0   \cr
 0  & \tau_1 &\ddots & \vdots \cr
 \vdots &\ddots & \ddots & 0 \cr
      0  & \cdots & 0      & \tau_m   \end{matrix} \right) 
; \tau_j \in \mathrm {SO}_\mathbb C(2)\right\},$$ 
$$  R(H_1) = \{ \mathrm{Id}\} \mathrm{ \  and \ }
 R(H_2) = \{\mathrm{Id},
 \epsilon\} \mathrm{ \ with \ } \epsilon =  \left (  \begin{matrix} 
  -1        &         &        &   0                   \cr
            &    0     &   1    &                       \cr
            &    1    &   0     &                       \cr
    0       &         &        &  \mathrm{Id}_{2m-3}   \end{matrix} \right). 
\qquad \quad \ \    \quad$$ 

\medskip
\noindent
Observe that in the case $H_2$, where $\epsilon$ is present, $S$ is not
connected. But the ``$\epsilon$-part'' of $S$ is not relevant, 
since any $h\in H_2$ does not change the sign
of the first component of the imaginary part of $z_j \in T$ and therefore 
$\Sigma_{H_2\epsilon}(z)$ is empty for all $z\in T^N$. 
Thus it is sufficient to prove the 
following 

\begin{prop} For every possible $H\in \{H_0,H_1,H_2\}$ and 
every $w \in T^N$ the set $\Sigma_H(w) = \{ h \in H ;
h \cdot w \in T^N \}$ is connected.
\end{prop}

\medskip
\noindent
\textit{Proof.} We will carry out the proof in the case where 
$n=2m-1$ and $H=H_0$. The proof in the other cases is analogous. 
Note that $H$ splits into its real and imaginary part, 
i.e., $H = H_\mathbb R \cdot H_I \cong H_\mathbb R \times H_I$ where 
$H_\mathbb R$ denotes the connected component  of the identity of 
$\lo 
\cap H = \{ h \in H ; 
\bar h = h \} $ and  $ H_I = \exp i \mathfrak h_\mathbb R$.
Thus the $2\times 2$ blocks appearing for  $h\in H_I$ are given by 
$$ \sigma = \left( \begin{matrix} a  & i b \cr
i b & a \end{matrix} \right) \ \text{where} \  a^2+b^2 =1 
\mathrm{ \quad and \quad}
\tau_j =   \left( \begin{matrix} c_j & -i d_j  \cr
i d_j  & c_j  \end{matrix} \right) \ \text{where}\ c_j^2 - d_j^2 = 1, c_j > 0.
$$
Let $S^1:=\{(x,y)\in\RR^2;\ x^2+y^2=1\}$,  
$\mathcal H:=\{(x,y) \in \mathbb R^2;\ x^2 - y^2 = 1 
\mbox{ and } x>0\}$, identify $H_I$ with $S^1  \times \mathcal H  \times 
\cdots \times \mathcal H\subset \RR^2\times \cdots\times \RR^2=\RR^{2m}$ 
and let $$\tilde \psi:\RR^{2m}\to \RR^{(1+n)\times (1+n)},\ 
\tilde \psi(a,b,c_1, d_1, \ldots, c_{m-1}, d_{m-1})=\left( \begin{matrix} 
\sigma &  0 & \cdots & 0 \cr
     0      & \tau_1 & \ddots & \vdots \cr
 \vdots & \ddots & \ddots & 0 \cr
       0    & \cdots & 0 &   \tau_{m-1}   \end{matrix} \right)$$
where $\sigma = \left( \begin{matrix} a  & i b \cr
i b & a \end{matrix} \right)$ and 
$\tau_j =   \left( \begin{matrix} c_j & -i d_j  \cr
i d_j  & c_j  \end{matrix} \right)$. The restriction $\psi$ of $\tilde\psi$
to $S^1  \times \mathcal H  \times \cdots \times \mathcal H$ is a 
diffeomorphism onto its image $H_I$.

\medskip
\noindent
For every $w_k\in T$, $k=1, \ldots, N$  we get the  linear map
$ \tilde\varphi_k :\RR^{2m}\to\RR^{1+n}$, $p\mapsto 
\mathrm{Im}(\tilde \psi(p)\cdot w_k)$. Note that
\begin{itemize}
\item  If
$p= (p_1,..,p_m)\in\tilde \varphi_k^{-1}(C)$, then 
$(p_1,..,rp_j,..,p_m)\in\tilde \varphi_k^{-1}(C)$ for all $0 < r \leq 1$
and $j=2,..,m$.

\item  If $p = (p_1,..,p_m), p_j \in\tilde \varphi_k^{-1}(C)$, then 
$(s\cdot p_1,p_2,..,p_m)\in  \tilde \varphi_k^{-1}(C)$ for all $s > 1$.
\end{itemize}
where $p_1=(a,b), p_j=(c_j,d_j)\in\RR^2, j=2,\ldots, m$.

\medskip
\noindent
It remains to show that $\Sigma_{H_I}(w)$ is connected for all 
$w \in T^N$. 

\medskip
\noindent
Let  $e := ((1,0),(1,0),...,(1,0) )= \psi^{-1}(\mathrm{Id})
 \in \psi^{-1}(\Sigma_{H_I}(w))$ and $p = (p_1,..,p_m) := \psi^{-1}(h)
\in \psi^{-1}(\Sigma_{H_I}(w))$. From the convexity of $C$ and the 
linearity of $\tilde \varphi_k$ it  follows that $ q(t) = (q_1(t),..,q_m(t))=
e + t(p-e)$ is contained in $\bigcap_{k= 1}^N \varphi_k^{-1} (C)$ 
for $t \in [0,1]$.
Thus 
$$ \tilde  \gamma_p(t) := \left ( \frac {q_1(t)} {\Vert q_1(t) \Vert_E},
\frac {q_2(t)} { \sqrt{ \eta(q_2(t))}}, ..., 
\frac { q_m(t)} {\sqrt{\eta (q_m(t))}} \right) \in \psi^{-1} ( \Sigma_H(w)) $$ 
for
$t \in [0,1]$. Here $\Vert \cdot  \Vert_E$ denotes the standard Euclidean norm.
 Thus   $\gamma_h(t) := \psi(\tilde \gamma_p(t))$ gives a curve 
which connects Id with $h$.
\qd

\medskip
\noindent
Since $\lo$ is a real form of $\loc$, orbit connectness implies the following 
(see \cite{H1})

\begin{corr} \label{equivariant}
Let $Y$ be a complex space with a holomorphic $\loc$-action.
Then every holomorphic $\lo$-equivariant  map $ \varphi : T^N \to Y$ extends 
to a holomorphic $\loc$-equivariant map $\Phi : \loc \cdot T^N \to Y$. \qd
\end{corr}

\medskip
\noindent
In the terminology of \cite{H1} Corollary \ref{equivariant}
 means that $\loc \cdot T^N$ is
the universal complexification of the $\lo$-space $T^N$.

\section{The strictly plurisubharmonic exhaustion of the tube} \label{spsh}

\medskip
\noindent
Let $X,Q,P$ be topological spaces, $q : X \to Q$ and $ p : X \to P$
continuous maps. A function $f : X \to \RR$ is said to be an exhaustion of 
$X$ mod $p$ along $q$ if for every compact subset $K$ of $Q$ and $r \in \RR$
the set $p(q^{-1}(K) \cap f^{-1}((-\infty,r]))$ is compact.

\medskip
\noindent
The characteristic function of the forward cone $C$ is up to a
constant given by the function
$ \tilde \rho : C \to \RR, \tilde \rho (y) = \eta(y)^{ -\frac {n+1} 2}$.
It follows from the construction of the characteristic function,
that $\log \tilde \rho $ is a $\lo$-invariant strictly convex 
function on $C$ (see \cite{Far} for details).
In particular
$$ \rho : T^N \to \mathbb R , \quad (x_1+iy_1,..,x_N+iy_N) \mapsto \frac 
1 {\eta (y_1)} + \dots + \frac 1 { \eta (y_N) }$$
is a $\lo$-invariant strictly plurisubharmonic function on $T^N$.
Of course this may also be checked by direct computation.

\medskip
\noindent Let $\pi_\CC : \vn \to \vn \mo \loc$ be the analytic Hilbert 
quotient and $\pi_\RR : T^N \to T^N / \lo$ the quotient by the $\lo$-action.
In the following we always write $z = x+iy$, i.e., $z_j = x_j +iy_j$ where
$x_j$ denote the real and $y_j$ the imaginary part of $z_j$. For example 
$ z_j \bullet z_k = x_j \bullet x_k - y_j \bullet y_k + i ( x_j \bullet y_k 
+ x_k \bullet y_j).$

\medskip
\noindent
The main result of this section is the following 

\begin{theo} \label{exhaustion}
The function $\rho : T^N \to \RR$,
is an exhaustion of $T^N$ mod $\pi_\RR$ along $\pi_\CC$.
\end{theo}

\medskip
\noindent
We do the case of one copy first.

\begin{lemm} \label {1copy}
 Let $D_1 \subset T$ and assume that $\pi_\CC (D_1) \subset \CC$ is
 bounded. Then  $\{ ( x\bullet y,$  $\eta(x), \eta(y) ) \in \RR^3 
; z =  x+iy \in D_1\}$
  is bounded.
\end{lemm}

\medskip
\noindent
\textit{Proof.}
 The condition on $D_1$ means, that there is a $M \geq 0$ 
such that 
$$ \vert \eta(x) -  \eta(y) \vert \leq M \qquad \mbox {and} \qquad 
\vert x \bullet y \vert \leq M$$
for all $z = x+iy \in D_1$. Since $\eta(x) \eta(y) \leq (x \bullet y)^2$ and
$\eta(y)  \geq 0$, this implies that 
$\{ (x \bullet y, \eta(x) ,\eta(y)) \in \RR^3 ; z \in D_1\}$ is bounded.
\qd

\begin{lemm}
 Let $D_2 \subset T \times T$ be such that $\pi_\CC(D_2)$ is
bounded. Then $\{ (\eta(x_1),\eta(y_1),$ $\eta(x_2),$ $\eta(y_2),$ $x_1 
\bullet x_2$, $y_1 \bullet y_2) \in \RR^6 ; (z_1,z_2) \in D_2\}$ is bounded.
\end{lemm}

\medskip
\noindent
\textit{Proof.}
 Lemma \ref{1copy} implies that there is a $M_1 \geq 0$ such that 
$\vert \eta(x_j)\vert \leq M_1, \vert \eta(y_j) \vert \leq M_1$ and 
$\vert x_j \bullet y_j\vert \leq M_1$,  $j = 1,2$, for all ($z_1,z_2) \in D_2$.
Now $\eta(z_1+z_2) = \eta (z_1) + \eta (z_2) + 2 \cdot z_1 \bullet z_2$
shows that $\{ \eta(z_1 +z_2) \in \RR ; (z_1,z_2) \in D_2\}$ is bounded.
 But $z_1 +z_2 \in T$, thus  Lemma \ref{1copy}
 implies $\vert \eta (x_1+x_2) \vert
\leq M_2$ and $ \vert \eta(y_1+y_2) \vert \leq M_2$ for some $M_2 \geq 0$
  and all $(z_1,z_2) \in D_2$. This gives 
$$ \vert x_1 \bullet x_2 \vert \leq \frac 32 \mbox{ max }\{M_1,M_2\} 
\qquad \mbox{and} \qquad
 \vert y_1 \bullet y_2 \vert \leq \frac 32 \mbox{ max }\{M_1,M_2\}.$$
\qd

\begin{beme}  \label{mixed terms}
Based on the following we  only need, that the set
$\{(\eta(y_1), \eta(y_2),
y_1 \bullet y_2) \in \RR^3;$ $(z_1,z_2) \in D_2\}$ is bounded. 
We apply this to 
points $y_j + iy_1$ where 
$\pi_\CC(y_j + iy_1) = \eta(y_j) - \eta(y_1) +2i y_j \bullet y_1.$
\end{beme}

\begin{beme} \label{Slice}
For every subset $X$ of \ $ T,$ we have 
$$ X \subset \lo \cdot (X \cap ( \RR^{1+n} +i ( \RR^{>0} \cdot e_0))),$$
where  $\RR^{>0} \cdot e_0 = \{ t e_0; t > 0 \} \subset \RR^{1+n}$.
\end{beme}

\begin{lemm} \label{M(B,K)}
 For every compact sets $B \subset C$ and $K \subset \CC$
the set $$M(B,K) := \{ x \in \RR^{1+n} ; \pi_\CC (x+iy) \in K 
\mbox{  for some } y  \in B\}$$
 is compact.
\end{lemm}

\medskip
\noindent
\textit{Proof.}
 Since $B$ and $K$ are compact, 
$M(B,K)$ is closed. We have to 
show that it is bounded. First note that $B_1 \subset B_2$ implies 
$M(B_1,K) \subset M(B_2,K)$. Using the properness of the 
$\lo$-action on $C$, we see, that there is an interval
$I = \{ t \cdot e_0 ; a \leq t \leq b \}$, $a > 0$ in $\RR \cdot e_0$ and
a compact subset $N$ in $\lo$, such that $N \cdot I := \bigcup_{g \in N}
 g \cdot I \supset B$. Thus $M(B,K) \subset M(N\cdot I,K) = N \cdot 
M(I,K) :=  \bigcup_{g \in N} g \cdot M(I,K)$. 

\medskip
\noindent
It remains to show that 
$M(I,K)$ is bounded. For $x \in M(I,K), x = \left( \begin{matrix}  x_0 
\cr \vdots \cr x_n \end{matrix} \right)$, there exists a $M_1 \geq 0$ such 
that $ \vert x \bullet 
(y_0 \cdot e_0)\vert = \vert x_0 \cdot y_0 \vert \leq M_1$ for all $y_0 \cdot 
e_0 \in I$. Since 
$a \leq y_0 \leq b$ and $a > 0$, this implies $\vert x_0^2 \vert \leq \frac 
{{M_1}^2}{ \vert y_0^2 \vert} \leq \frac {{M_1}^2} {a^2}$.
There also exists a $M_2 \geq 0$ such that 
$\vert \eta(x) \vert = \vert x_0^2 -x_1^2- \cdots x_n^2\vert \leq M_2$, so 
 we get $x_1^2+ \cdots x_n^2 \leq \frac {M_1^2} {a^2} + M_2$.
\qd

\begin{corr} \label{corollary}
For every $ r >0$ the set 
$M(B,K) \cap \{ y \in \RR^{1+n} ;  r \leq \eta(y) \}$ 
is compact. \qd
\end{corr}

\medskip
\noindent
\textit{Proof of  Theorem \ref{exhaustion}.} 
Using Remark \ref{Slice}  it is sufficient to prove that the set 
$$S := (\pi_\CC^{-1} (K) \cap \{ \rho \leq r\}) \cap ((\RR^{1+n} + i(\RR^{>0} 
\cdot e_0))  \times T^{N-1})$$ is compact. For $z = (z_1,..,z_N) \in S$ let
$z_j = x_j + iy_j$, where $x_j$ denotes the real part and $y_j$ the imaginary 
part of $z_j$. By the definition of $S$ we have $y_1 = y_{10} \bullet e_0$
where $y_{10}$ = $y_1 \cdot e_0$.  Moreover,
we get $\frac 1 r \leq \eta(y_1) = (y_{10})^2 \leq M$. Therefore the set 
$\{ y_1 \in \RR^{1+n}; (z_1,..,z_N) \in S\} 
= \{ t \cdot e_0 ; t^2 \in [\frac 1 r, M], t > 0\}$
is compact. 

\medskip
\noindent 
By Remark \ref{mixed terms} we get  that the sets 
$\{ (\eta(y_1),\eta(y_j), y_1 \bullet y_j)\in \RR^3; (z_1,..,z_N) \in S\}$
are bounded for $j=2,..,N$. Therefore we get the boundedness of
$\{ \pi_\CC (y_j +i y_1) \in  \CC;$ $(z_1,..,z_N) \in S \}$. Thus
the  $y_j$, $j = 2,..,N$,
with $(z_1,..,z_N) \in S$ are lying in the sets $M(I,B_j) \cap
  \{ y \in \RR^{1+n} ;  r \leq \eta(y) \}$, where 
$I := \{ t \cdot e_0 ; t^2 \in [\frac 1 r , M], t > 0 \}$ and $B_j$ are compact 
subsets of $\CC$, containing $\{ \pi_\CC (y_j +i y_1) \in  \CC;$ $(z_1,..,z_N)
 \in S \}$. By Corollary  \ref{corollary} these sets are compact, which 
implies that the set 
$\{ (y_1,..,y_N) \in \RR^{(1+n) \times N} ; (z_1,..,z_N) \in S\}$ 
is compact. Hence using  Lemma  \ref{M(B,K)} it follows that 
$\{ (x_1,..,x_N) \in \RR^{(1+n) \times N}; (z_1,..,z_N) \in S\}$
is bounded. Thus $S$ is bounded and therefore compact.
\qd

\section{  Saturatedness of the extended future tube} \label{sat}

\medskip
\noindent
We call $A \subset X$ saturated with respect to a map $p : X \to Y$ if 
$A$ is the inverse image of 
a subset of $Y$. 

\medskip
\noindent
Let $\pi_\mathbb C : \vn \to \vn \mo \loc$  be the analytic Hilbert quotient, 
which is given by the algebra of $\loc$-invariant polynomials
functions on $\vn$ (see section \ref{GIT}) and let $U_r $ denote the set $\{
z \in T^N ; \rho(z) < r\}$ for some  $r  \in \RR \cup \{+\infty\},$ 
where $\rho$ is the 
strictly plurisubharmonic exhaustion function, which we defined
in  section \ref{spsh}.

\begin{theo} \label{saturated}
The set  $ \loc \cdot  U_r = \loc \cdot  \{ z \in T^N ; \rho(z) < r\}$ is 
saturated with respect to $\pi_\mathbb C$.
\end{theo}

\medskip
\noindent
It  is well known, that  each fiber of $\pi_\mathbb C$ contains exactly 
one closed orbit of $\loc$ (see section \ref{GIT}).  Moreover, every 
orbit contains a closed orbit 
in its closure. Therefore it is sufficient
to prove

\begin{prop}
If $z \in U_r$ and $\loc \cdot u$ is the  closed orbit in 
$\overline{\loc \cdot z}$, then $\loc \cdot u \cap U_r \neq \emptyset$. 
\end{prop}

\medskip
\noindent
The idea of proof is to construct a  one-parameter group $\gamma$
 of $\loc$, such that $\gamma(t) z \in U_r $ for $ \vert t \vert \leq 1$ and 
$\lim_{t \to 0} \gamma(t) z \in \loc \cdot u$.

\medskip
\noindent
In the following, let $z=(z_1,..,z_N) \in U_r$ and denote by  
$L(z) = \CC z_1 + \cdots +  \CC z_N $ the $\CC$-linear subspace 
of $\mathbb C^{1+n}$ spanned by $z_1,..,z_N$. The subspace of isotropic 
vectors in $L(z)$ with respect to the Lorentz product is denoted by $L(z)^0$,
i.e., $L(z)^0= \{ w \in L(z); w \bullet v = 0$ for all $ v \in L(z) \}$. 
Let  $\overline{L(z)^0}$ be its  conjugate, i.e., 
$\overline{L(z)^0} = \{ \bar v ; v \in L(z)^0\}$.   

\begin{lemm}
For all  $\omega  \neq 0$, $\omega \in L(z)^0$ we have 
$\eta(\mathrm{Im}(\omega)) < 0$.
\end{lemm}

\medskip
\noindent
\textit{Proof.}
 Let $\omega = \omega_1 +i \omega_2$ with $\omega_1 = \mathrm{Re} 
(\omega), \omega_2 = \mathrm{Im} (\omega)$. Assume that 
 $\eta(\mathrm{Im}(\omega)) =  \eta(\omega_2) \geq 0$. Since $\omega 
\in L(z)^0$, we have
 $0 = \eta(\omega) = 
\eta(\omega_1) - \eta(\omega_2) + 2i \omega_1 \bullet \omega_2$.

\medskip
\noindent
If $ \eta(\omega_2) > 0$, i.e., $\omega_2  \in C$ or $\omega_2  \in -C$,
then $\omega_1 \bullet \omega_2 = 0$ contradicts $\eta(\omega_1) = 
\eta (\omega_2) >0$. Thus assume $\eta(\omega_1) = \eta (\omega_2) =0$
and $\omega_1 \bullet \omega_2 = 0$. Hence $\omega_1$ and
$\omega_2$ are $\RR$-linearly dependent and therefore there is a 
$\lambda \in \CC, \omega_3 \in \RR^{1+n}$ such that 
$\omega = \lambda \omega_3$ and $\omega_3 \bullet e_0 \geq 0$. 
We have $\eta(\omega_3) = 0$ and, since $\omega_3 \in L(z)^0, e_0 \bullet 
\omega_3 \geq 0$ and $z_1 \in T$, we also have $ 0 = \omega_3 \bullet 
\mathrm{Im} (z_1)$. This implies by   the definition of $C$ that $\omega_3 =0$.
\qd

\begin{corr}
For $\omega \in L(z)^0, \omega \neq 0$, we have 
$\omega \bullet \bar \omega <0$. In particular, 
$L(z)^0 \cap \overline{ L(z)^0} = \{0\}$ and  the complex Lorentz
product is non-degenerate on $L(z)^0 \oplus  \overline{ L(z)^0}.$ \qd
\end{corr}

\begin{corr} \label{W}
Let $W := (L(z) \oplus \overline {L(z)})^\perp :=
 \{ v \in \CC^{1+n}; v \bullet u = 0$ for all $ u \in 
L(z)^0 \oplus  \overline{ L(z)^0}\}$. Then 
$$ L(z) = L(z)^0 \oplus ( L(z) \cap W).$$ \qd
\end{corr}

\medskip
\noindent
\textit{Proof of Proposition \ref{saturated}.}
Let $z  \in U_r$. We use the notation of Corollary \ref{W}.
Define 
$$\gamma : \mathbb C^\ast \to \loc \mathrm{\quad by 
\quad}
\gamma(t) v = \begin{cases} t v & \mbox{ for }  v \in L(z)^0\cr 
		       t^{-1} v & \mbox{ for }  v \in \overline{L(z)^0} \cr
	     v & \mbox{ for }  v \in W
  \end{cases}. $$
Every component $z_j$ of $z$ is of the form 
$z_j = u_j + \omega_j$ where  $u_j \in W$ and $ \omega_j \in L(z)^0$
are uniquely determined by $z_j$. Recall that $W$ is the set  
$\{ v \in \CC^{1+n}; v \bullet u = 0$ for all $ u \in 
L(z)^0 \oplus  \overline{ L(z)^0}\}$. Since $\lim_{t \to 0} \gamma(t) z_j
= u_j$ and $L(u)^0 = \{0\}$ for $ u = (u_1,..,u_N)$, $u$ lies in the unique 
closed orbit in $\overline{ \loc.z}$ (see Lemma \ref{closedorbits}).
It remains to show that $u \in U_r$. For every $t \in \CC$
 we have 
$$ \eta(\mathrm{ Im }(u_j + t\omega_j)) = 
\eta (\mathrm{ Im } (u_j))  + 
 \vert t \vert^2 \eta(\mathrm{Im}(\omega_j)).$$ Since 
$\eta (\mathrm{ Im } (u_j +  \omega_j))>0$ and 
$\eta(\mathrm{Im}(\omega_j)) \leq 0$, this implies 
$\eta(\mathrm{Im}(u_j + t \omega_j)) \in C^\pm$ for all $ t  \in 
[0,1].$ Moreover, 
$ \eta(\mathrm {Im}(z_j)) < \eta(\mathrm{Im}(u_j)),$ for every $j$. Thus
$\rho(z) > \rho(u)$ and therefore $u \in U_r.$
\qd

\begin{corr} 
The extended future tube is saturated with respect to $\pi_\CC$. \qd
\end{corr}

\begin{beme} \label{konvex}
The function 
$f : \RR \to \RR, t \mapsto \eta(\mathrm{ Im} (u_j + t \omega_j))$, 
is strictly concave if $\omega_j \neq 0$. The proof shows
 $u_j + t \omega_j \in T$
 for all $t \in \RR.$
\end{beme}

\section{The K\"ahlerian reduction of the extended future tube} 
\label{Kreduction}

\medskip
\noindent
If one is only interested in the statement of the future tube conjecture,
one can simply apply the main result in \cite{H2} (Theorem 1 in $\S$2). 
Our goal here is to  show that much more is true.

\medskip
\noindent
For $z \in \vn$ let $x = \frac 12 (z + \bar z)$ be the real and
$y = \frac 1{2i}(z - \bar z)$ the imaginary part of $z$, i.e., $z=(z_1,..,z_N)
= (x_1,..,x_N)+i(y_1,..,y_N)$ in the obvious sense.
The strictly plurisubharmonic function $\rho : T^N \to \RR,$
$\rho(z) = \frac 1 {\eta(y_1)} + \cdots + \frac  1 {\eta(y_N)}$
defines for every $\xi \in \mathfrak{so}(1,n) = \mathfrak o (1,n)$
the function $$\mu_\xi(z) = d\rho(z) (i\xi z) =  \frac  {d}
{dt} \Big \vert_{t=0} \  \rho(\exp it \xi \cdot z).$$
Here of course $\mathfrak {so}(1,n) = \mathfrak o(1,n)$ denotes the Lie 
algebra of $\mathrm O_\RR(1,n)$. The real group $\lo$ acts by conjugation on 
$\mathfrak {so}(1,n)$ and therefore by duality on the dual vector space 
$\mathfrak {so}(1,n)^\ast$. It is easy to check that the map 
$\xi \to \mu_\xi$ depends linearly on $\xi$. Thus $$\mu : T^N \to 
\mathfrak {so}(1,n)^\ast, \quad \mu(z)(\xi) := \mu_\xi(z),$$ is 
a well defined $\lo$-equivariant map. In fact $\mu$ is a moment map
with respect to the K\"ahler form $\omega = 
2i \partial  \bar \partial  \rho.$ 

\medskip
\noindent
In order to emphasizes the general ideas, we set $G := \lo$,
$\gc := \loc$, $X:= T^N$ and $Z:= \gc \cdot X$. The corresponding analytic
Hilbert quotient, induced by $\pi_\CC : \vn \to \vn \mo \loc$ are denoted
by $\pi_X : X \to X \mo G$, $\pi_Z : Z  \to Z \mo \gc$. Note that,
by what we proved, we have $X \mo G = Z \mo \gc$.

\begin{prop}  \ 
\begin{enumerate} 
\item For every $q \in Z \mo \gc $ we have $(\pi_\CC)^{-1}(q) \cap \mu^{-1}(0)
= G \cdot x_0$ for some $x_0 \in \mu^{-1}(0)$ and $\gc \cdot x_0$ 
is a closed orbit in $Z.$

\item The inclusion $\mu^{-1}(0) \stackrel \iota  \to X \subset Z$ induces a 
homeomorphism $\mu^{-1}(0)/G  \stackrel {\bar \iota}  \to Z \mo \gc.$
\end{enumerate}
\end{prop}

\medskip
\noindent
\textit{Proof.} A simple calculation shows that the set of critical points of 
$\rho \vert \gc \cdot x \cap X,$ i.e., $\mu^{-1}(0) \cap \gc \cdot x,$ consists of
a discrete set of $G$-orbits. Moreover, every critical point is a local 
minimum (see \cite{H2}, Proof of Lemma 2 in $\S$ 2).

\medskip
\noindent
On the other hand Remark \ref{konvex} of section \ref{sat} says that if 
$\rho \vert \gc \cdot x \cap X$ has a local minimum in $x_0 \in \gc \cdot 
x \cap X,$
then $\gc \cdot x_0 = \gc \cdot x$ is necessarily closed in $Z$. Moreover,
$\rho \vert \gc \cdot x \cap X$ is then an exhaustion and therefore 
$\mu^{-1}(0) \cap ( \gc \cdot x_0 \cap X) = G \cdot x_0$ 
(see \cite{H2}, Lemma 2 in $\S$ 2).
This proves the first part.

\medskip
\noindent The statement i. implies that 
$\iota : \mu^{-1}(0) \hookrightarrow X \subset Z$ 
induces a bijective continuous map $\bar \iota : \kn /G \to Z \mo  \gc.$ 
Since the $G$-action on $X$ is proper and $\kn$ is closed,
the action on $\kn$ is proper. In particular $\kn/G$ is
a Hausdorff topological space.

\medskip
\noindent 
Theorem \ref{saturated} implies that $\bar \iota$ is a homeomorphism, since for 
every sequence $q_\alpha \to q_0$ in $Z\mo\gc$ we find a sequence 
$(x_\alpha)$ such that $x_\alpha$ are contained in a 
compact subset of $\kn$ and $\pi_\CC(x_\alpha) = q_\alpha.$ Thus 
every convergent subsequence of $(x_\alpha)$ has  a limit point 
in $G \cdot x_0$ where  $\pi_\CC(x_0) = q_0$.
\qd

\begin{prop} \label{rho}
The restriction $\rho \vert \kn : \kn \to \RR$ induces a strictly 
plurisubharmonic continuous exhaustion $\bar \rho : Z \mo \gc \to \RR. $
\end{prop}

\medskip
\noindent
\textit{Proof.} The exhaustion property for $\bar \rho$ follows from Theorem
\ref{exhaustion}. The argument that $\bar \rho$ is strictly
plurisubharmonic is the same as in \cite{HHL}.
\qd

\begin{theo}
The extended future tube $Z$ is a domain of holomorphy.
\end{theo}

\medskip
\noindent
\textit{Proof.}
 Proposition \ref{rho} implies that $Z \mo \gc$ is a Stein space (see  
\cite{N} Theorem II).
Hence $Z$ is a Stein space.
\qd

\medskip
\noindent
In fact, much more has been proved here. We would like to comment on this.
By definition, an analytic subset of a complex manifold is closed.
For the following recall that orbit-connectedness is a condition
on the $\gc$-orbits.

\begin{prop}
Every analytic $G$-invariant subset 
$A$ of $X$ is orbit connected in $Z$ and $\gc \cdot A$ is 
an analytic subset of $Z$. In particular, $\gc \cdot A$ is a Stein space.
Moreover the restriction maps
$$ \mathcal O(Z)^{\gc} \to \mathcal O (\gc \cdot A)^{\gc} \to \mathcal O(A)^G$$
are surjective.
\end{prop}

\medskip
\noindent
\textit{Proof.}
If $b \in \gc \cdot A \cap X$, then $ b = g \cdot a$ 
for some $g \in \gc$ and $a \in A$. Hence $g \in \Sigma_{\gc} (a)
= \{ g \in \gc ; g \cdot a \in X\}$. The 
identity principle for holomorphic functions shows that
$ \Sigma_{\gc}(a) \cdot a \in A$. Thus $b \in A$ This shows 
$ \gc \cdot A \cap X = A$. But $\{ g \cdot  X ; g \in \gc\}$ is an open 
covering of $X$ such that $\gc \cdot  A \cap g \cdot  X = g \cdot  A$. 
This shows that $\gc \cdot A$ is 
an analytic subset of $Z$. In particular, it is a Stein space. 
The last statement follows from orbit connectedness (see \cite{H1}).
\qd

\begin{prop}
For every $G$-invariant analytic subset $A$, its saturation
$\hat A = \pi_X^{-1}(\pi_X(A))$ 
is an analytic subset of $X$. Moreover, $\hat A \mo G$ is canonically 
isomorphic to $A \mo G$ and 
$\pi_{\hat A}: \hat A \to \hat A \mo G \subset X \mo G$ 
is the Hilbert quotient of $\hat A$ whose restriction to $A$ 
gives the analytic Hilbert quotient of $A$
\end{prop}

\medskip
\noindent
\textit{Proof.}
We already know that $A^c = \gc \cdot A$ is an 
analytic subset of $Z$. Its saturation $\hat A^c = \pi_Z^{-1} ( \pi_Z(A^c))
=  \pi_Z^{-1} ( \pi_Z(A))$ is
an analytic subset of $Z$ and it is easily checked that 
$\hat A = \hat A^c \cap X = \pi_X^{-1} ( \pi_X(A))$ has the desired
properties.
\qd

\bigskip
\noindent
P. Heinzner

\noindent
Fakult\"at und Institut  f\"ur Mathematik

\noindent
Ruhr-Universit\"at Bochum

\noindent
Geb\"aude NA 4/74

\noindent
D-44780 Bochum

\noindent
Germany

\noindent
e-mail: heinzner@cplx.ruhr-uni-bochum.de

\bigskip
\noindent
P. Sch\"utzdeller

\noindent
Fakult\"at und Institut f\"ur Mathematik 

\noindent
Ruhr-Universit\"at Bochum

\noindent 
Geb\"aude NA 4/69

\noindent
D-44780 Bochum

\noindent
Germany

\noindent
e-mail: patrick@cplx.ruhr-uni-bochum.de

\begin{thebibliography}{9999999}


\bibitem[B]{Brem} {\bf R. Bremigan},
{\rm Invariant analytic domains in complex semisimple groups},
{\it Transformation Groups} {\bf 1} (1996), 279--305

\bibitem[FK]{Far} {\bf J. Faraut, A. Koranyi},
{\rm Analysis on Symmetric Cones},
Oxford Press, Oxford  1994


\bibitem[HW]{HW} {\bf D. Hall, A. D. Wightman},
{\rm A theorem on invariant analytic functions with applications
to relativistic quantum field theory},
{\it Kgl. Danske Videnskap. Selkap, Mat.-Fys. Medd} {\bf 31} (1965) 1--14

\bibitem[He1]{H1} {\bf P. Heinzner},
{\rm Geometric invariant theory on Stein spaces},
{\it Math. Ann.} {\bf 289} No. 4 (1991), 631--662

\bibitem[He2]{H2} {\bf P. Heinzner},
{\rm The minimum principle from a Hamiltonian point of view},
{\it Doc. Math. J.} {\bf 3} (1998), 1--14

\bibitem[HeHuL]{HHL} {\bf P. Heinzner, A. T. Huckleberry, F. Loose},
{\rm K\"ahlerian extensions of the symplectic reduction},
{\it J. reine angew. Math.} {\bf 455} (1994), 123--140

\bibitem[HeMP]{HMP} {\bf P. Heinzner, L. Migliorini, M. Polito},
{\rm Semistable quotients},
{\it Ann. Scuola Norm. Sup. Pisa} {\bf 26} (1998), 233--248

\bibitem[J]{Jo} {\bf R. Jost},
{\rm The general theory of quantized fields},
In: Lectures in applied mathematics vol. IV,  1965



\bibitem[Kr]{Kr} {\bf H. Kraft},
{\rm Geometrische Methoden in der Invariantentheorie},
In: Aspects of Mathematics, Vieweg Verlag 1984


\bibitem[N]{N} {\bf R. Narasimhan},
{\rm The Levi Problem for Complex Spaces II},
{\it Math. Ann.} {\bf 146} (1962), 195--216

\bibitem[SV]{SV} {\bf A. G. Sergeev, V. S. Vladimirov},
{\rm Complex analysis in the future tube},
In: Encyclopaedia of mathematical sciences (Several complex variables
 II) vol. 8 (1994), 179--253

\bibitem[StW]{SW} {\bf  R. F. Streater, A. S.  Wightman},
{\rm PCT spin statistics, and all that},
W. A. Benjamin, INC. 1964


\bibitem[W]{W} {\bf A. S.  Wightman},
{\rm Quantum field theory and analytic functions of several complex variables},
{\it J. Indian Math. Soc.} {\bf  24} (1960), 625--677


\bibitem[Z]{Z} {\bf X. Y. Zhou},
{\rm A proof of the extended future tube conjecture},
{\it Izv. Math.} {\bf 62} (1998), 201--213

\end{thebibliography}
\end{document}